\documentclass{amsart}
\usepackage{amssymb,amsmath,amsthm,xspace}%\usepackage{showkeys}
\usepackage[small]{diagrams}
\def\mat#1{\ensuremath{#1}\xspace}

\def\makemath#1#2#3{%#1-new command, #2-new box, #3-math expression
\newsavebox{#2}
\sbox{#2}{\ensuremath{#3}}
\def#1{\usebox{#2}\xspace}}

\def\cA{\mat{\mathbb{A}}}   %affine space

\def\cF{\mat{\mathbb{F}}}
\def\cN{\mat{\mathbb{N}}}   %natural numbers
\def\cQ{\mat{\mathbb{Q}}}   %ring of rationals
\def\cC{\mat{\mathbb{C}}}   %complex numbers

\def\cZ{\mat{\mathbb{Z}}}   %ring of integers

\def\lE{\mat{\mathcal{E}}}

\def\lH{\mat{\mathcal{H}}}

\def\lL{\mat{\mathcal{L}}}
\def\lM{\mat{\mathcal{M}}}

\def\lP{\mat{\mathcal{P}}}

\def\lX{\mat{\mathcal{X}}}
\def\lY{\mat{\mathcal{Y}}}

\makemath{\Psi}{\Psibox}{\Psi}
\makemath{\Phi}{\Phibox}{\Phi}

\def\la{\mat{\lambda}}
\def\vi{\mat{\varphi}}
\def\si{\mat{\sigma}}

\def\om{\mat{\omega}}

\def\al{\mat{\alpha}}

\def\ga{\mat{\gamma}}
\def\Ga{\mat{\Gamma}}

\def\hi{\mat{\chi}}

\def\ksi{\mat{\xi}}

\def\te{\mat{\theta}}

\def\mrm@#1{\mat{\mathrm{#1}}}

\def\g {\mat{\mathfrak{g}}}  %\gg and \g are standartly reserved for ">>" and ">>>"
\def\DMO{\DeclareMathOperator}

\DMO{\Hom}{Hom}
\DMO{\lHom}{\lH\mathit{om}}
\DMO{\Ext}{Ext}
\DMO{\lExt}{\lE\mathit{xt}}
\DMO{\End}{End}
\DMO{\Aut}{Aut}
\DMO{\Fun}{Fun}
\DMO{\Tor}{Tor}
\DMO{\ext}{ext}

\DMO{\Ob}{Ob}
\DMO{\Mor}{Mor}
\DMO{\im}{im}
\DMO{\coim}{coim}
\DMO{\coker}{coker}
\DMO{\Arr}{Arr}

\DMO{\Id}{Id}
\DMO{\add}{add} % splitting of idempotents (karoubinization)
\DMO{\ind}{ind} % category of ind-objects
\DMO{\pro}{pro} % category of pro-objects
\DMO{\Map}{Map}
\DMO{\Iso}{Iso}
\DMO{\Isom}{Isom}

\DMO{\Presh}{Presh}

\DMO\coalg{Coalg}
\DMO{\Rep}{Rep}

\DMO{\Mod}{Mod}
\DMO{\rad}{rad}
\DMO{\soc}{soc}
\DMO{\ann}{ann}

\DMO{\Spec}{Spec}
\DMO{\spec}{Spec}
\DMO{\Proj}{Proj}
\DMO{\supp}{supp}
\DMO{\Coh}{Coh}
\DMO{\coh}{Coh}
\DMO{\Qcoh}{QCoh}
\DMO{\QCoh}{QCoh}
\DMO{\Pic}{Pic}
\DMO{\Div}{Div}
\DMO{\ch}{ch}
\DMO{\Hilb}{Hilb}
\DMO{\Fitt}{Fitt}
\DMO{\Quot}{Quot}
\def\Gm{\mat{{{\mathbb G}_{\mathrm m}}}}
\DMO{\Gras}{Gr}
\DMO{\Flag}{Flag}

\DMO{\cone}{cone}
\DMO{\Tw}{Tw}
\DMO{\rank}{rk}
\DMO{\rk}{rk}
\DMO{\codim}{codim}
\DMO{\cov}{cov}
\DMO{\sgn}{sgn}
\DMO{\td}{td}
\DMO{\GL}{GL}
\DMO{\SL}{SL}
\def\gl{\mat{\mathfrak{gl}}}

\DMO\Der{Der}
\DMO\der{Der}
\DMO\coder{Coder}
\DMO{\diag}{diag}
\DMO{\HMod}{HMod} %the homotopy category of modules over DGC
\DMO{\ad}{ad}
\DMO*{\colim}{colim}
\DMO*{\hocolim}{hocolim}
\DMO*{\holim}{holim}
\DMO{\Ho}{Ho}
\DMO{\har}{char}
\DMO{\sk}{sk}
\DMO{\cosk}{cosk}
\DMO{\Gal}{Gal}
\DMO{\tr}{tr}
\DMO{\Tr}{Tr}
\DMO{\Sh}{Sh}
\DMO{\Is}{Is} %Isometries
\DMO{\Hol}{Hol} %Holomorphic automorphisms
\DMO{\Lie}{Lie} %Lie algebra of a group
\DMO{\Res}{Res} %restriction
\DMO{\irr}{irr} %
\DMO{\Irr}{Irr} %
\DMO{\Exp}{Exp} %
\DMO{\Log}{Log} %
\DMO{\mult}{mult} %
\DMO{\height}{ht}

\newcommand{\GIT}{/\!\!/}

\def\dd{\mat{\partial}}

\def\iso{\simeq}

\def\sb{\subset}

\def\xx{\times}

\def\lpoly#1{[\![#1]\!]} %formal power series [[#1]]

\def\inv{^{-1}}

\def\ub#1{\mat{\overline{#1}}}  %upper bar
\def\set#1{\mat{\{ #1\}}}
\def\sets#1#2{\mat{\{ #1 \mid #2\}}}

\def\arrowsUsual{
\newarrow{TeXto}----{->}
\newarrow{TeXinto}C---{->}
\newarrow{TeXonto}----{->>}
\def\ar{\rightarrow}
\def\emb{\hookrightarrow}
\def\mto{\mapsto}
\def\arr{rTeXto}
\def\embb{\rTeXinto}
\newarrow{Eq}=====
    }

\newif\ifukr\ukrfalse
\newif\ifrus\rusfalse

\def\theorems{
\newtheorem{prp}{\ifukr Пропозиція \else \ifrus Предложение \else Proposition\fi\fi}[section]
\newtheorem{proposition}[prp]{\ifukr Пропозиція \else \ifrus Предложение \else Proposition\fi\fi}

\newtheorem{theorem}[prp]{\ifukr Теорема \else \ifrus Теорема \else Theorem\fi\fi}

\newtheorem{lemma}[prp]{\ifukr Лема \else \ifrus Лемма \else Lemma\fi\fi}

\newtheorem{corollary}[prp]{\ifukr Висновок \else \ifrus Следствие \else Corollary\fi\fi}

\theoremstyle{definition}

\newtheorem{definition}[prp]{\ifukr Означення \else \ifrus Определение \else Definition\fi\fi}

\newtheorem{remark}[prp]{\ifukr Зауваження \else \ifrus Замечание \else Remark\fi\fi}
}

\theorems\arrowsUsual

\def\field{\mat\cF}
\def\Rep{\mathcal R}
\def\xrt{\al} %root vector
\def\xwt{\mat{\nu}} %weight vector

\usepackage{ifpdf}
\ifpdf
\usepackage[pdftex]{hyperref}
\else
\usepackage{hyperref}
\fi

\begin{document}
\title[]{On the multiplicities of the irreducible
highest weight modules over Kac-Moody algebras}%

\author{Sergey Mozgovoy}%

\address{Institut f\"ur Mathematik, Johannes Gutenberg-Universit\"at Mainz,
55099 Mainz, Germany.}%

\email{mozgov@mathematik.uni-mainz.de}%

\thanks{}%

%\date{10.09.2006}%
%\dedicatory{}%
%\commby{}%
\maketitle
\begin{abstract} We prove that the weight multiplicities of the
integrable irreducible highest weight module over the Kac-Moody
algebra associated to a quiver are equal to the root
multiplicities of the Kac-Moody algebra associated to some
enlarged quiver. To do this, we use the Kac conjecture for
indivisible roots and a relation between the Poincar\'e
polynomials of quiver varieties and the Kac polynomials, counting
the number of absolutely irreducible representations of the quiver
over finite fields. As a corollary of this relation, we get an
explicit formula for the Poincar\'e polynomials of quiver
varieties, which is equivalent to the formula of Hausel
\cite{Haus1}.
\end{abstract}

\renewcommand{\theprp}{\arabic{prp}}
\section{Introduction}
Let $(\Ga,I)$ be a finite quiver without loops, where $I$ is the
set of vertices. The underlying graph of $\Ga$ defines a symmetric
generalized Cartan matrix $C$ with $c_{ij}$ being equal to minus
the number of edges connecting the vertices $i$ and $j$ if $i\ne
j$ and $c_{ii}=2$. Let $\g(\Ga)$ be the corresponding Kac-Moody
algebra \cite{Kac3}. We will identify $\cZ^I$ with its root
lattice. Define the quadratic form $T$ on $\cZ^I$, called the Tits
form, by the matrix $\frac 12C$. For any $\xwt=(\xwt_i)_{i\in
I}\in\cN^I$, we define $\ub\xwt:=\sum_i\xwt_i\om_i$, where
$(\om_i)_{i\in I}$ are the fundamental weights (which are in
general not unique). Our purpose is to relate the multiplicities
of the integrable irreducible highest weight module $L(\ub\xwt)$
with the root multiplicities of the Kac-Moody algebra associated
to some enlarged quiver.

We define this enlarged quiver $\Ga_*$ by adjoining to $\Ga$ a new
vertex $*$ and $\xwt_i$ arrows from $*$ to $i$, for each $i\in I$.
Its set of vertices is $I_*=I\cup\set{*}$. As before, we may
consider the Kac-Moody algebra $\g(\Ga_*)$ associated to the
quiver $\Ga_*$. For any $\xrt\in\cZ^I$ and any $k\in\cZ$, we
consider the pair $(\xrt,k)$ as an element of $\cZ^{I_*}$. Our
main result is the following

\begin{theorem}\label{thr intro multiplicities} For any
$\xrt\in\cN^I$, we have $\dim
L(\ub\xwt)_{\ub\xwt-\xrt}=\dim\g(\Ga_*)_{(\xrt,1)}$.
\end{theorem}

For example, it allows to perform an explicit computation of the
multiplicities of $L(\ub\xwt)$ using the Peterson recursive
formula for the root multiplicities, see \cite{Pet1,M2}. It also
allows to write down the character formula for the part
$\oplus_\xrt \g(\Ga_*)_{(\xrt,1)}$ of the Kac-Moody algebra
$\g(\Ga_*)$ and one can conjecture that similar formulas for the
other parts of $\g(\Ga_*)$ exist as well.

To prove the theorem, we use quiver varieties. Let $T_*$ be the
Tits form of the quiver $\Ga_*$ and let
$d=d(\xrt,\xwt):=1-T_*(\xrt,1)=\xrt\cdot\xwt-T(\xrt)$, where
$\xrt\cdot\xwt=\sum_i\xrt_i\xwt_i$. It is known that $\dim L(\ub
\xwt)_{\ub\xwt-\xrt}$ equals $h_c^{2d}(\lM)$, where
$\lM=\lM(\xrt,\xwt)$ is a quiver variety having dimension $2d$,
see \cite{Nak2,Saito1}. Using the trick of Crawley-Boevey
\cite{CB1} one can consider the quiver variety $\lM(\xrt,\xwt)$ as
a certain moduli space of representations of the double quiver
$\ub\Ga_*$ of the quiver $\Ga_*$ (obtained by adjoining reverse
arrows for all arrows in $\Ga_*$). This allows us to relate the
Poincar\'e polynomial of $\lM$ (and thus $h^{2d}_c(\lM)$) to
certain invariants of the quiver $\Ga_*$.

Namely, for any quiver $(\Ga,I)$ and any $\xrt\in\cN^I$, let
$a_\xrt(\Ga,q)$ be the number of absolutely indecomposable
representations of $\Ga$ over $\cF_q$ of dimension \xrt. It is
proved in \cite{Kac1} that $a_\xrt(\Ga,q)$ are polynomials in $q$
with integer coefficients. Moreover, $a_\xrt(\Ga,q)\ne0$ if and
only if \xrt is a root of $\g(\Ga)$ and $a_\xrt(\Ga,q)=1$ if and
only if $\xrt$ is a real root. We call these polynomials the Kac
polynomials of the quiver \Ga. It was conjectured by Kac
\cite{Kac1} that $a_\xrt(\Ga,0)=\dim\g(\Ga)_\xrt$. This conjecture
was proved for indivisible roots in \cite{CBB}. The proof of the
full conjecture was announced by Hausel \cite{Haus1}. We derive
Theorem \ref{thr intro multiplicities} from the Kac conjecture for
indivisible roots and our second result, which is

\begin{theorem}\label{thr intro rel between Poinc and Kac poly}
We have $\sum_i h^{2i}_c(\lM)q^i=q^{d} a_{(\xrt,1)}(\Ga_*,q)$ and
$h^{2i+1}_c(\lM)=0$ for $i\ge0$.
\end{theorem}

The proof of the theorem is based on the paper \cite{CBB}. In the
course of the proof we additionally show that the polynomial
$p(\lM,q):=\sum_i h^{2i}_c(\lM)q^i$ counts the number of points in
$\lM(\cF_q)$ for a finite field $\cF_q$ of a sufficiently large
characteristic.

Notice, that the Poincar\'e polynomial (of cohomologies with
compact support) of \lM equals $p(\lM,t^2)$. Theorem \ref{thr
intro rel between Poinc and Kac poly} allows to calculate the
Poincar\'e polynomial of \lM using a formula for $a_\xrt(\Ga,q)$
for arbitrary quiver $(\Ga,I)$, see \cite[Theorem 4.6]{Hua1} or
\cite[Theorem 2]{M2}. Following \cite{SB1}, we define the function
$r_\xrt(\Ga,q)$ by the formula
$$\frac 1{\#\GL_\xrt(\cF_q)}%
\#\sets{(g,x)\in(\GL_\xrt\xx\Rep(\Ga,\xrt))(\cF_q)}{gx=x,\ g\text{
is unipotent}},$$%
where $\GL_\xrt$ and $\Rep(\Ga,\xrt)$ are defined in Section
\ref{sec moduli of repres}. There is an explicit formula for
$r_\xrt(\Ga,q)$ and an easy relation between the generating
functions
$$a(\Ga,q):=\sum_{\xrt\in\cN^I}a_\xrt(\Ga,q)x^\xrt\qquad\text{and}\qquad%
r(\Ga,q):=\sum_{\xrt\in\cN^I}r_\xrt(\Ga,q)x^\xrt,$$%
see Section \ref{sec Poincare poly}. Using the relation between
$a(\Ga,q)$ and $r(\Ga,q)$ we obtain from Theorem~\ref{thr intro
rel between Poinc and Kac poly}

\begin{theorem}\label{thr intro hausel formula} We have
$$\sum_{\xrt}q^{-d(\xrt,\xwt)}p({\lM(\xrt,\xwt)},q)x^\xrt%
=(q-1)\frac{\sum_\xrt r_{(\xrt,1)}(\Ga_*,q)x^{\xrt}}{r(\Ga,q)}.$$
\end{theorem}

This formula is equivalent to the formula of Hausel \cite[Theorem
5]{Haus1}. As well as Theorem \ref{thr intro rel between Poinc and
Kac poly} it gives an effective way to calculate the Poincar\'e
polynomials of quiver varieties.

In the second section we recall the definition of the moduli
spaces of quiver representations according to \cite{King1}. In the
third section we recall the basic properties of quiver varieties.
Section \ref{sec main results} is devoted to the proof of Theorem
\ref{thr intro multiplicities} and Theorem \ref{thr intro rel
between Poinc and Kac poly}. In Section \ref{sec Poincare poly} we
recall the relation between Kac polynomials and function
$r(\Ga,q)$ and then prove Theorem \ref{thr intro hausel formula}.

\renewcommand{\theprp}{\thesection.\arabic{prp}}

\section{Moduli spaces of quiver representations}\label{sec moduli of repres}
In this section we follow closely \cite{CBB}. Let $\Ga$ be a
quiver and $I$ be its set of vertices. For any arrow $h\in \Ga$,
we denote by $h'$ and $h''$ its source and target respectively. We
denote by $\ub\Ga$ the double of $\Ga$, obtained from it by
adjoining reverse arrows for all arrows in $\Ga$. For any $h\in
\ub\Ga$, we denote by $\ub h$ the opposite arrow from $\ub\Ga$.
For any $\xrt,\xwt\in\cZ^I$ we define $\xrt\cdot \xwt:=\sum_{i\in
I}\xrt_i\xwt_i$.

Let $R$ be a commutative ring. For any $I$-graded free $R$-module
$V$ of finite rank, we define $\dim V:=(\rk V_i)_{i\in
I}\in\cN^I$. Given two $I$-graded free $R$-modules $V$ and $W$, we
denote the module of $I$-graded morphisms between them by
$\Hom_I(V,W)$.

Let $\xrt\in\cN^I$ and let $V$ be an $I$-graded free $R$-module
with $\dim V=\xrt$. Define the scheme over $R$
$$\Rep(\Ga,\xrt):=\bigoplus_{h\in\Ga}\Hom(V_{h'},V_{h''}).$$
Then we can identify $\Rep(\ub \Ga,\xrt)$ with
$\Rep(\Ga,\xrt)\oplus \Rep(\Ga,\xrt)^*$. There is an obvious
action of the group
$$\GL_\xrt:=\prod_{i\in I}\GL_{\xrt_i}$$
on $\Rep(\Ga,\xrt)$ and therefore on $\Rep(\ub \Ga,\xrt)$ (this
action can be factored through $G_\xrt=\GL_\xrt\!\!/\Gm$, where
$\Gm$ is considered as a diagonal subgroup in $\GL_\xrt$). There
is a map (which is a moment map if $R$ is a field)
$$\mu:\Rep(\ub \Ga,\xrt)\ar \g_\xrt^*\emb\gl_\xrt^*$$
defined by
$$(x_h)_{h\in\ub\Ga}\mto \sum_{h\in \Ga}[x_h,x_{\ub h}],$$
where $\gl_\xrt=\prod_i M_{\xrt_i\xx\xrt_i}$ is a Lie algebra of
$\GL_\xrt$ and $\gl_\xrt^*$ is isomorphic to $\gl_\xrt$ by the
trace pairing; $\g_\xrt$ is a Lie algebra of $G_\xrt$ and
$\g_\xrt^*\emb\gl_\xrt^*$ can be identified with such matrices
$(\ksi_i)_{i\in I}$ that $\sum_i\tr\ksi_i=0$.

\begin{definition} We call an element $x\in\Rep(\Ga,\xrt)$
nilpotent, if there exists some $N\ge1$ such that for any path
$h_1\dots h_N$ in \Ga (i.e. $h''_i=h'_{i+1}$, $1\le i<N$) it holds
$x_{h_N}\dots x_{h_1}=0$.
\end{definition}

Now we recall some facts from the paper of King \cite{King1} about
the moduli spaces of semistable representations of quivers. In his
paper King uses the geometric invariant theory over an
algebraically closed field. According to Seshadri \cite{Seshadri1}
this can be done also over \cZ. The quotients obtained by Seshadri
are categorical quotients but in general not universal categorical
quotients.
%(and this is needed in order to be sure that after the
%base change the quotients are obtained by the same procedure as
%before).
This problem was overcome in \cite[Lemma B.4]{CBB}, where
it was shown that after base change from \cZ to some $\cZ_f$,
$f\in\cZ$ the quotients are universal categorical quotients. This
implies that we can consider these quotients over \cC and over
finite fields of a sufficiently large characteristic. In what
follows, we will avoid all this formalities and describe the
constructions over an algebraically closed field, but we will bear
in mind that everything can be done over $\cZ_f$. This will be
used later.

Let \field be an algebraically closed field. For any $\te\in\cZ^I$
with $\te\cdot \xrt=0$, we define a character
$\hi_\te:\GL_\xrt\ar\Gm$ by the formula $\hi_\te(g)=\prod_{i\in
I}\det(g_i)$. This character defines an action of $\GL_\xrt$ on
the trivial line bundle $L$ over $\Rep(\Ga,\xrt)$.

\begin{definition} A point $x\in \Rep(\Ga,\xrt)$ is called
\te-stable (respectively, \te-se\-mi\-stable) if for any
$I$-graded, $x$-invariant subspace $0\ne V'\subsetneq V$ it holds
$\te\cdot\dim V'>0$ (respectively, $\te\cdot \dim V'\ge0$).
\end{definition}

It is proved in \cite{King1} that the stable part
$\Rep(\Ga,\xrt)^s$ (respectively, semistable part
$\Rep(\Ga,\xrt)^{ss}$) of $\Rep(\Ga,\xrt)$ with respect to the
$\GL_\xrt$-line bundle $L$ consists precisely of \te-stable
(respectively \te-semistable points). By the geometric invariant
theory (see \cite{King1} for details) there exists a geometric
quotient $\Rep(\Ga,\xrt)^s/G_\xrt$ and a categorical quotient
$\Rep(\Ga,\xrt)^{ss}\GIT G_\xrt$. Moreover, the inclusion
$\Rep(\Ga,\xrt)^{ss}\emb \Rep(\Ga,\xrt)$ induces a projective map
$\Rep(\Ga,\xrt)^{ss}\GIT G_\xrt\ar \Rep(\Ga,\xrt)\GIT G_\xrt$. In
the same way one defines stability conditions and the
corresponding moduli spaces for the representations of $\ub \Ga$.

%Here we use a slightly modified notion of semistability from
%\cite[Section 3]{CBB}. For any $\te\in\cZ^I$ (called further a
%character), we define the slope function
%$s(\xrt):=\frac{\te\cdot\xrt}{\height\xrt}$, where
%$\height\xrt=\sum\xrt_i$. Define a point $x\in R(\Ga,\xrt)$ to be
%\te-stable (\te-semistable) if for any $I$-graded, $x$-invariant
%subspace $0\ne V'\subsetneq V$ it holds $s(\dim V')<s(\dim V)$
%($s(\dim V')<s(\dim V)$). We denote the set of \te-stable points
%by $R(\Ga,\xrt)^s$ and denote the set of \te-semistable points by
%$R(\Ga,\xrt)^{ss}$. It is proved in \cite{King1} that there exists a
%geometric quotient $R(\Ga,\xrt)^s/G_\xrt$ and a categorical quotient
%$R(\Ga,\xrt)^{ss}\GIT G_\xrt$. Moreover, the inclusion
%$R(\Ga,\xrt)^{ss}\emb R(\Ga,\xrt)$ induces a projective map
%$R(\Ga,\xrt)^{ss}\GIT G_\xrt\ar R(\Ga,\xrt)\GIT G_\xrt$. In the same way
%one defines stability conditions for the representations of $\ub
%\Ga$.

\begin{lemma} The moment map $\mu:\Rep(\ub
\Ga,\xrt)^s\ar\g_\xrt^*$ is smooth.
\end{lemma}
\begin{proof} Stability condition implies that the stabilizer in
$G_\xrt$ of any stable point is trivial. This implies that $\mu$
is smooth at any stable point.
\end{proof}

\begin{corollary} If $\mu\inv(0)^s$ is nonempty then the map
$\mu:\Rep(\ub \Ga,\xrt)^s\ar\g_\xrt^*$ is surjective.
\end{corollary}
\begin{proof} It is known that the smooth morphisms are open, so
the image of the map $\mu:\Rep(\ub \Ga,\xrt)^s\ar\g_\xrt^*$ is
open. As this image contains $0$ and is stable with respect to the
multiplication by scalars, it coincides with $\g_\xrt^*$.
\end{proof}

\begin{lemma}\label{lem equality for fibers} Assume that
$\Rep(\ub \Ga,\xrt)^s=\Rep(\ub \Ga,\xrt)^{ss}$ and there exists
some $G_\xrt$-invariant element $\ksi\in\g_\xrt^*$ with
$\mu\inv(\ksi)^s\ne\emptyset$. Then the map $\mu:\Rep(\ub
\Ga,\xrt)^s\ar\g_\xrt^*$ is surjective and for finite fields
\field of a sufficiently large characteristic we have
$\#\mu\inv(\ksi)^s\GIT G_\xrt(\field)=\#\mu\inv(0)^s\GIT
G_\xrt(\field)$.
%the complex varieties $\mu\inv(\ksi)^s\GIT
%G_\xrt(\cC)$ and $\mu\inv(0)^s\GIT G_\xrt(\cC)$ are homotopic(No?)
\end{lemma}
\begin{proof} We have to prove just the second statement. It is a
consequence of \cite[Appendix by Nakajima]{CBB}. Let $L$ be a line
through $0$ and $\ksi$ in $\g_\xrt^*$, $X:=\mu\inv(L)^s$ and
$Y=\mu\inv(L)$. Let $\lX=X\GIT G_\xrt$ and $\lY:=Y\GIT G_\xrt$.
There is a commutative diagram \begin{diagram}
X&\rInto&Y\\
\dTo&&\dTo \\
\lX&\rTo^\pi&\lY&\rTo&L,
\end{diagram}
where $\pi$ is projective. Let $\Gm$ acts on $\Rep(\Ga,\xrt)$ by
multiplication of all the matrices by scalar. This action induces
the action of $\Gm$ on $X,\ Y,\ \lX,\ \lY$. There is also an
action of $\Gm$ on $L$ s.t. the maps in the above diagram are
$\Gm$-equivariant. We claim that for any $x\in \lX$ there exists
$\lim_{t\to 0}tx$. First of all, the map $\Gm\ar \lY$, $t\mto
t\pi(x)$ can be extended to $\cA^1$ by $0\mto 0$. Now it follows
from the projectivity of $\pi$ that the map $\Gm\ar\lX$, $t\mto
tx$ can also be extended to $\cA^1\ar\lX$. This proves the
existence of the limit. We apply now the result of \cite[Appendix
by Nakajima]{CBB} to the smooth family $\lX\ar L$ with the above
$\Gm$-action and get the second statement.
\end{proof}

\section{Quiver varieties}
We assume that \field is an algebraically closed field. Let
$\xwt\in \cN^I$ be fixed and let $W$ be an $I$-graded vector space
of dimension \xwt. Let $\Ga_*$ be the quiver defined in
Introduction; that is, we adjoin to $\Ga$ a new vertex $*$ and
$\xwt_i$ arrows from $*$ to $i$ for each $i\in I$. We identify
$W_i$ with $\bigoplus_{h:*\ar i}\field\cdot h$. Let as before $V$
be an $I$-graded vector space of dimension \xrt. As in
Introduction, a pair $(\xrt,n)$ with $\xrt\in\cN^I$ and $n\in\cN$
will be considered as an element from $\cN^{I_*}$.

There is an obvious identification
$$M(\xrt,\xwt):=\Rep(\ub \Ga_*,(\xrt,1))=\Rep(\ub \Ga,\xrt)\oplus\Hom_I(W,V)\oplus\Hom_I(V,W).$$
The elements of this space will be represented as triples
$(x,p,q)$. Note that $G_{(\xrt,1)}=(\prod_{i\in
I}\GL_{\xrt_i}\xx\Gm)/\Gm\iso\GL_\xrt$. Therefore the moment map
may be considered as a map $\mu_*:M(\xrt,\xwt)\ar\gl_\xrt^*$. It
is given by the formula
$$\mu_*(x,p,q)=\mu(x)+pq.$$

We fix once and for all $\te\in \cZ^{I_*}$,
$\te=(-1,\dots,-1,\sum\xrt_i)$ and consider stability and
semistability conditions in $M(\xrt,\xwt)$ with respect to \te.

\begin{lemma} Stability and semistability conditions in
$M(\xrt,\xwt)$ are equivalent. An element $(x,p,q)\in
M(\xrt,\xwt)$ is stable if and only if any $I$-graded,
$x$-invariant subspace $V'\sb V$ s.t. $q(V')=0$ is zero.
\end{lemma}
\begin{proof} Assume that $(x,p,q)$ is semistable and there
exists an $I$-graded, $x$-invariant subspace $V'\sb V$ s.t.
$q(V')=0$. Then $V'$ may be considered as an $I_*$-graded,
$(x,p,q)$-invariant subspace of $V\oplus\field$ of dimension
$(\dim V',0)$. From semistability condition we get $-\dim V'\ge0$,
thus $V'=0$. It follows that the last condition of the lemma
holds.

Conversely, assume that the last condition of the lemma holds. Let
$V'\oplus V_*$ be some proper, $I_*$-graded, $(x,p,q)$-invariant
subspace of $V\oplus\field$. If $V_*=\field$, then $\dim V'<\xrt$
and therefore $\te\cdot(\dim V',1)=\sum_i\xrt_i-\sum_i\dim
V'_i>0$. If $V_*=0$, then $q(V')=0$ and by our assumption $V'=0$.
This implies that $(x,p,q)$ is stable.
\end{proof}

\begin{remark} For any subscheme $X\sb M(\xrt,\xwt)$ we write
respectively $X^n$, $X^s$, $X^{ns}$ to denote the subschemes of
$X$ consisting respectively of nilpotent, stable, nilpotent and
stable elements.
\end{remark}

\begin{definition} We define the quiver variety
$\lM=\lM(\xrt,\xwt)$ to be the quotient $\mu_*\inv(0)^s\GIT
\GL_\xrt$. Define $\lL=\lL(\xrt,\xwt):=\mu_*\inv(0)^{ns}\GIT
\GL_\xrt$.
\end{definition}

\begin{remark} It is easy to see that $\lL(\xrt,\xwt)$ is the
preimage of zero under the projective morphism $\mu_*\inv(0)^s\GIT
\GL_\xrt\ar \mu_*\inv(0)\GIT \GL_\xrt$. It is known that an
element $(x,p,q)\in M(\xrt,\xwt)^s$ is nilpotent if and only if
$x$ is nilpotent and $p=0$, see e.g.\ \cite[Lemma 5.9]{Nak1} or
\cite[Lemma 2.22]{Lus1}.
\end{remark}

Let $T$ denotes the Tits form of the quiver \Ga and $T_*$ denotes
the Tits form of the quiver $\Ga_*$. As in Introduction, we define
$d=d(\xrt,\xwt):=1-T_*(\xrt,1)=\xrt\cdot\xwt-T(\xrt)$.

\begin{theorem}[{Nakajima \cite[Section 3]{Nak2}}] Variety $\lM$
is smooth and variety $\lL$ is projective. The complex manifold
$\lM(\cC)$ is symplectic and its subvariety $\lL(\cC)$ is a
Lagrangian subvariety homotopic to $\lM(\cC)$. The dimension of
$\lM$ equals $2d(\xrt,\xwt)$.
\end{theorem}

For a scheme $X$ of finite type over $\cZ_f$, we define
$$h^i(X):=\dim H^i(X(\cC),\cQ),\qquad h^i_c(X):=\dim H^i_c(X(\cC),\cQ).$$
Note that $h^i_c(X)$ can also be defined as $\dim
H^i_c(X_{\ub\cF_p},\cQ_l)$ (for large enough prime $p$) by the
base change theorem \cite[Thorem 1.8.7]{FK1} and comparison
theorem \cite[Thorem 1.11.6]{FK1}.

\begin{lemma}\label{lem mult of modules} $\dim
L(\ub\xwt)_{\ub\xwt-\xrt}=h^{2d}_c(\lM)$.
\end{lemma}
\begin{proof} It is well-known (see e.g.\ \cite{Nak2} or
\cite{Saito1}) that $\dim L(\ub\xwt)_{\ub\xwt-\xrt}$ equals the
number of irreducible components of \lL i.e., $h^{2d}_c(\lL)$. We
note that
$$h^{2d}_c(\lL)=h^{2d}(\lL)=h^{2d}(\lM)=h^{2d}_c(\lM),$$
where the last equality follows from the Poincar\'e duality.
\end{proof}

\section{Main results}\label{sec main results}
Recall from Introduction that for any $\xrt\in\cN^I$ there is a
polynomial $a_\xrt(\Ga)\in\cZ[q]$ such that for any finite field
$\cF_q$, $a_\xrt(\Ga,q)$ is the number of absolutely
indecomposable representations of $\Ga$ over $\cF_q$ of dimension
\xrt. As before $\lM=\lM(\xrt,\xwt)$ is a quiver variety.

\begin{proposition}\label{prp number of points of qv} For a
finite field $\cF_q$ of a sufficiently large characteristic it
holds $\#\lM(\cF_q)=q^{d(\xrt,\xwt)}a_{(\xrt,1)}(\Ga_*,q)$.
\end{proposition}
\begin{proof} We are going to apply \cite[Proposition
2.2.1]{CBB}. Let $\xrt':=(\xrt,1)\in\cN^{I_*}$. Note that
$\te=(-1,\dots,-1,\sum\xrt_i)$ is $\xrt'$-generic in the sense
that for any $0<\ga'<\xrt'$ it holds $\te\cdot \ga'\ne 0$. Indeed,
let $\ga'=(\ga,k)$, where $\ga\in\cZ^I$ and $k\in\set{0,1}$. If
$k=0$ then clearly $\te\cdot\ga'=-\sum_{i}\ga_i<0$. If $k=1$ then
$\ga<\xrt$ and $\te\cdot\ga'=\sum_i\xrt_i-\sum_i\ga_i>0$. We
identify \te with an element of $\gl_{\xrt'}$ consisting of
diagonal matrices. Then all the points of $\mu_*\inv(\te)$ are
stable. Indeed, assume that some $(x,p,q)$ is not stable i.e.,
there exists an $I$-graded, $x$-invariant subspace $0\ne
V'\subsetneq V$ with $q(V')=0$. Then it follows that $\sum_{h\in
\Ga}[x_h,x_{\ub h}]\big|_{V'}=(-\Id_{V'_i})_{i\in I}$. Adding the
traces we get $0=\sum_i\dim V'_i$, which is impossible. Applying
now \cite[Proposition 2.2.1]{CBB}, we get
$a_{\xrt'}(\Ga_*,q)=q^{-d(\xrt,\xwt)}\cdot\#\mu_*\inv(\te)\GIT\GL_\xrt(\cF_q)$
for a finite field $\cF_q$ of a sufficiently large characteristic.
By Lemma \ref{lem equality for fibers} we have
$\#\mu_*\inv(\te)\GIT\GL_\xrt(\cF_q)=\#\lM(\cF_q)$.
\end{proof}

It follows that there exists a polynomial $p(\lM)\in\cZ[q]$ such
that for finite fields $\cF_q$ of large enough characteristic it
holds $\#\lM(\cF_q)=p(\lM,q)$.

\begin{proof}[Proof of Theorem \ref{thr intro rel between Poinc
and Kac poly}] According to \cite[Lemma A.1]{CBB} and
\cite[Proposition A.2]{CBB} it holds
$p(\lM,q)=\sum_{i}h_c^{2i}(\lM)q^i$ and $h^{2i+1}_c(\lM)=0$
whenever we show the existence of a \Gm-action on \lM s.t. for any
$x\in \lM$ there exists $\lim_{t\to 0}tx$ and s.t. $\lM^{\Gm}$ is
projective. This action was described in Lemma \ref{lem equality
for fibers}. Note that the \Gm-invariant part of \lM is mapped to
zero under the map $\pi:\lM\ar\mu_*\inv(0)\GIT\GL_\xrt$. This
implies that $\lM^\Gm\sb\lL$ and therefore $\lM^\Gm$ is
projective.
\end{proof}

\begin{proof}[Proof of Theorem \ref{thr intro multiplicities}] By
the Kac conjecture, proved for indivisible vectors in \cite{CBB},
it holds $\dim\g(\Ga_*)_{(\xrt,1)}=a_{(\xrt,1)}(\Ga_*,0)$. By
Proposition \ref{prp number of points of qv} it holds
$a_{(\xrt,1)}(\Ga_*,0)=q^{-d(\xrt,\xwt)}p(\lM,q)\big|_{q=0}$. So,
we have to prove $\dim
L(\ub\xwt)_{\ub\xwt-\xrt}=q^{-d(\xrt,\xwt)}p(\lM,q)\big|_{q=0}$.
From the facts that $\lM$ is homotopic to $\lL$ and that \lL is
projective of dimension $d(\xrt,\xwt)$ we get that
$h^i(\lM)=h^i(\lL)=h^i_c(\lL)=0$ for $i>2d(\xrt,\xwt)$ and
therefore $h^i_c(\lM)=0$ for $i<2d(\xrt,\xwt)$ by Poincar\'e
duality. It follows that
$h^{2d(\xrt,\xwt)}_c(\lM)=q^{-d(\xrt,\xwt)}p(\lM,q)\big|_{q=0}$
and we apply Lemma \ref{lem mult of modules}.
\end{proof}

\section{Poincar\'e polynomials of quiver varieties}\label{sec Poincare poly}
In this section we recall the explicit formula for the functions
$r_\xrt(\Ga,q)$ defined in Introduction and the relation between
$r(\Ga,q)$ and $a(\Ga,q)$, see \cite{Hua1,M2}. From this relation
and Theorem \ref{thr intro rel between Poinc and Kac poly} we
derive then Theorem \ref{thr intro hausel formula}.

Let $(\Ga,I)$ be a finite quiver. Let \lP be the set of partitions
(see e.g.\ \cite{Mac1}) and let $\la=(\la^i)_{i\in I}\in\lP^I$.
Define $|\la|:=(|\la^i|)_{i\in I}\in\cN^I$. For any $j\ge1$ define
$\la_j:=(\la_j^i)_{i\in I}\in\cN^I$. Define $T(\la):=\sum_{j\ge
1}T(\la_j)$, where the quadratic form $T$ on $\cZ^I$ is the Tits
form defined in Introduction. Then the function
$r_\xrt(\Ga)\in\cQ(q)$, $\xrt\in\cN^I$ defined in Introduction
equals
$$r_\xrt(\Ga,q):=\sum_{|\la|=\xrt}\frac{q^{-T(\la)}}{\prod_{i\in I}\vi_{\la^i}(q\inv)},$$
where $\vi_\mu(q):=\prod_{j\ge 1}\vi_{\mu_j-\mu_{j+1}}(q)$ for
$\mu\in\lP$ and $\vi_n(q):=(1-q)\dots (1-q^n)$ for $n\in\cN$.

To describe the relation between
$a(\Ga)=\sum_{\xrt}a_\xrt(\Ga)x^\xrt$ and
$r(\Ga)=\sum_{\xrt}r_\xrt(\Ga)x^\xrt$, we use \la-rings (see e.g.\
\cite[Appendix]{M2}). We endow the field $\cQ(q)$ with the
structure of a \la-ring in terms of Adams operations by
$$\psi_n(f(q)):=f(q^n),\qquad f\in\cQ(q).$$
In order to avoid any problems with the Adams operations in what
follows, we tensor our \la-rings with \cQ without mentioning that
and so we always assume that our \la-rings contain \cQ. If $R$ is
a \la-ring, we endow the ring $R[x_1,\dots,x_r]$ with a \la-ring
structure by the formula $\psi_n(ax^\xrt):=\psi_n(a)x^{n\xrt}$,
where $a\in R$, $\xrt\in\cN^r$. In the same way, we endow with a
\la-ring structure the ring of formal power series over $R$. Given
a \la-ring $R$, we define the map $\Exp:R\lpoly
{x_1,\dots,x_r}^+\ar 1+R\lpoly {x_1,\dots,x_r}^+$ (here $R\lpoly
{x_1,\dots,x_r}^+$ is an ideal $(x_1,\dots,x_r)$) by the formula
$\Exp(f):=\sum_{n\ge0}\si_n(f)$ or, in terms of Adams operations,
$$\Exp(f)=\exp\big(\sum_{k\ge 1}\frac 1k\psi_k(f)\big),$$
where the map $\exp$ (as well as the map $\log$ used below) is
defined as in \cite[Ch.II \S 6]{BourLie13}. The map $\Exp$ has an
inverse $\Log:1+R\lpoly {x_1,\dots,x_r}^+\ar R\lpoly
{x_1,\dots,x_r}^+$ which is given by the formula of Cadogan (see
\cite{Cadogan1,Getz1,M2})
$$\Log(f)=\sum_{k\ge1}\frac {\mu(k)}k\psi_k(\log(f)),$$
where $\mu$ is a M\"obius function. Now we are ready to write down
the relation between $a(\Ga,q)$ and $r(\Ga,q)$ (see \cite[Theorem
2]{M2})
$$a(\Ga,q)=(q-1)\Log(r(\Ga,q)).$$
This formula together with Theorem \ref{thr intro rel between
Poinc and Kac poly} enables us to calculate the Poincar\'e
polynomial of quiver varieties. The content of Theorem \ref{thr
intro hausel formula} is a formula for the direct computation of
the Poincar\'e polynomial of quiver varieties.
\begin{proof}[Proof of Theorem \ref{thr intro hausel formula}]
For any $n\in\cN$, we define
$$a_n:=\sum_{\xrt\in\cN^I}a_{(\xrt,n)}(\Ga_*)x^\xrt,\qquad %
r_n:=\sum_{\xrt\in\cN^I}r_{(\xrt,n)}(\Ga_*)x^\xrt,$$
$a_*:=\sum_{n\ge0}a_nx_*^n$ and $r_*:=\sum_{n\ge0}r_nx_*^n$. Then
by Theorem \ref{thr intro rel between Poinc and Kac poly} it holds
$$\sum_{\xrt}q^{-d(\xrt,\xwt)}p(\lM(\xrt,\xwt),q)x^\xrt=a_1(q).$$
We know that $a_*(q)=(q-1)\Log(r_*(q))$ and therefore
\begin{multline*} a_1(q)=\frac \dd{\dd x_*}a_*(q)\big|_{x_*=0}
=(q-1)\frac \dd{\dd
x_*}\sum_{k\ge1}\frac {\mu(k)}k\psi_k(\log(r_*(q)))\big|_{x_*=0}\\
=(q-1)\frac \dd{\dd x_*}\log(r_*(q))\big|_{x_*=0}
=(q-1)\frac{\frac \dd{\dd x_*}r_*(q)}{r_*(q)}\big|_{x_*=0}
=(q-1)\frac {r_1(q)}{r_0(q)}.
\end{multline*}
We note that $r_0(q)=r(\Ga,q)$.
\end{proof}

%\bibliography{../tex/fullbib}
%\bibliographystyle{../tex/hamsplain}
\bibliography{fullbib}
%\bibliographystyle{hamsplain}
%GATHER{../tex/fullbib.bib}
\end{document}